\documentclass{amsart}
\usepackage[dvips]{graphicx}
\usepackage{amssymb}
\usepackage{latexsym}
\usepackage{xypic}

\newtheorem{theorem}{Theorem}
\newtheorem{lemma}{lemma}[section]

\theoremstyle{definition}
\newtheorem{definition}[lemma]{Definition}
\newtheorem{example}[lemma]{Example}
\newtheorem{examples}[lemma]{Examples}

\newtheorem{corollary}[lemma]{Corollary}
\newtheorem{proposition}[lemma]{Proposition}

\theoremstyle{remark}
\newtheorem{remark}[theorem]{Remark}

\numberwithin{equation}{section}

\begin{document}

\title{  \sc{Schottky groups cannot  act on} $\mathbb{P}^{2n}_{\mathbb{C}}$ {\sc as subgroups of} $PSL(2n+1,\Bbb{C})$ }
\author{Angel Cano}
\address{ Instituto de Matem\'aticas Unidad Cuernavaca, Universidad Nacional Autonoma de M\'exico, Avenida  Universidad sin n\'umero, Colonia Lomas de Chamilpa, Cuernavaca,  Morelos, M\'exico}
\email{angel@matcuer.unam.mx}


\subjclass{Primary 37F99, 32Q, 32M Secondary 30F40, 20H10, 57M60, 53C}



\maketitle

\begin{abstract}
In this paper we look at a special type of discrete subgroups of
$PSL_{n+1}(\Bbb{C})$ called Schottky groups. We develop some basic
properties of these  groups
 and their  limit set when $n > 1$, and we prove that  Schottky groups
only occur in odd dimensions, {\it i.e.}, they
cannot be realized as subgroups of $PSL_{2n+1}(\Bbb{C})$.
\end{abstract}

\section{Introduction}
\label{intro}
Schottky groups play a significant role in the theory of classical  Kleinian groups and
Riemann surfaces (see for instance \cite{marden, maskit1, maskit}). Their  analogues in
 higher dimensions were introduced   by  Nori \cite{nori} and Seade-Verjovsky
 \cite{sva}, though these groups were also known to N. Hitchin (see the commentary of  Nori in
 in \cite {nori}). These are a special type of discrete groups of automorphisms of
 complex projective spaces
having  non-empty  region of discontinuity, where the action is ``free''  with
compact quotient. Hence they are a rich source  for  complex  compact manifolds
equipped  canonically with a projective structure.  Schottky groups  also have
very interesting dynamics in their limit set, the complement of the  region of
discontinuity. Moreover, these groups are neither Fuchsian ({\it i.e.}, subgroups of
$PU(n,1)$) nor affine in general. Thus, if  we want to study Kleinian actions on higher
dimensional complex projective  spaces, Schottky groups provide a very nice starting point.

So far Schottky groups have been studied only  for odd-dimensional
projective spaces (in  \cite { nori, sva}). It is thus  natural to
ask whether Schottky groups exists in even dimensions. In this
paper we prove they do not:  Schottky groups cannot act by complex
automorphism  on $\Bbb{P}^{2n}_{\Bbb{C}}$.  Hence, in order to
construct discrete  groups of automorphisms of
$\Bbb{P}^{2n}_{\Bbb{C}}$ with a
 rich underlying geometry and  dynamics one must follow different methods. This is
 done
  for $\Bbb{P}^2_{\Bbb{C}}$ in  \cite {angel, pablo1, pablo2}.

This paper is divided into four sections. In section 1 we define
what  Schottky groups are and we state the main result of this
article.  In section 2  we develop some basic dynamical and
algebraic facts about Schottky groups. In section 3 we look at the
 limit set of infinite cyclic groups;    and in section 4  we use the previous information
 to show that Schottky groups cannot be realized  in even dimensions.

\section{Notations and the Main Result}
\label{sec:1}

We recall that the complex projective space
$\Bbb{P}^n_{\Bbb{C}}$ is defined as:
$$ \Bbb{P}^{n}_{\Bbb{C}}=(\Bbb{C}^{n+1}- \{0\})/\sim \,,$$
where "$\sim$" denotes the equivalence relation given by $x\sim y$
if and only if $x=\alpha y$ for some  non-zero complex scalar
$\alpha$. We know that $\Bbb{P}^n_{\Bbb{C}}$  is   a  compact
connected complex $n$-dimensional  manifold, which is naturally
equipped with the Fubini-Study  metric (see for instance
\cite{pablo2}).

If $[\mbox{ }]_n:\Bbb{C}^{n+1}-\{0\}\rightarrow \Bbb{P}^{n}_{\Bbb{C}}$ represents  the
quotient map, then a non-empty set  $H\subset \Bbb{P}^n_{\Bbb{C}}$ is said to be a
projective subspace of dimension $k$  (in symbols $dim_{\Bbb{C}} (H)=k$) if there is a
$\Bbb{C}$-linear subspace  $\tilde H$ of  dimension $k+1$
$($in symbols $dim_{\Bbb{C}}(\tilde H)=k+1)$, such that $[\tilde H]_n=H$.
Given a set of points $P$   in $\Bbb{P}^{n}_{\Bbb{C}}$, we define
$$\langle P \rangle =\bigcap\{l\subset \Bbb{P}^n_{\Bbb{C}}\mid  l \textrm{ is a
projective subspace and } P\subset  l \}.$$
So that  $\langle P\rangle $ is a projective subspace of $\Bbb{P}^{n}_{\Bbb{C}}$,
see \cite{manin}.\\

From now on, the symbols  $e_1,\ldots,e_{n+1}$ will either denote
the elements of the standard basis in $\Bbb{C}^{n+1}$ or  their
images under $[\mbox{ }]_n$.

Consider the action of $\Bbb{Z}_{n}$ (regarded as the $n$-roots of
unity) on $SL(n,\Bbb{C})$ given by $\alpha(a_{i,j})=(\alpha
a_{i,j})$. The quotient
$PSL_{n}(\Bbb{C})=SL_{n}(\Bbb{C})/\Bbb{Z}_{n}$ is a Lie Group
whose elements are called projective transformations. Every
representative   $\tilde \gamma$ of the coset
$\gamma=\Bbb{Z}_n\tilde \gamma=\gamma\in PSL_n(\Bbb{C})$ will be
called a lifting of $\gamma $. Observe that   $\gamma\in
PSL_{n+1}(\Bbb{C})$ acts  on $\Bbb{P}^n_{\Bbb{C}}$ as a
biholomorphic map by $\gamma([w]_n)=[\tilde\gamma(w)]_n$, where
$[w]_n\in \Bbb{P}^{n}_{\Bbb{C}}$ and $\tilde \gamma$ is a lifting
of $\gamma$.

\begin{definition} \label{d:s}
A subgroup $\Gamma\leq  PSL_{n+1}(\Bbb{C})$  is called a Schottky group if:

\begin{enumerate}
\item There are  $2g$ , $g\geq 2$, opens sets   $R_1,\ldots,R_g$, $S_1,\ldots,S_g$  in
$\Bbb{P}^n_{\Bbb{C}}$ with the  property  that:
\begin{enumerate}
\item each of these open sets  is the interior of  its closure; and
\item the closures  of the $2g$  open sets  are pairwise disjoint.
\end{enumerate}
\vskip.1cm
\item  $\Gamma$  has a generating set $Gen(\Gamma)=\{\gamma_1,\ldots,\gamma_g\}$ such that
  $\gamma(R_j)=\Bbb{P}^n_{\Bbb{C}}-\overline{S_j}$ for all $1 \leq j \leq g$, here the
  bar means topological closure.
\end{enumerate}
\end{definition}

From now on  $Int(A)$  will denote the topological  interior and
$\partial(A)$ the topological boundary of the set  $A$ and  for
each $1\leq j\leq g$, $R_j$ and $S_j$ will be denoted by
$R^*_{\gamma_j}$ and $S^*_{\gamma_j}$ respectively.

\begin{examples} \label{e:grou}
\begin{enumerate}
\item Every classical Schottky group of M\"obius  transformations
(see \cite{marden, maskit1, maskit}) is  Schottky in the sense of
definition \ref{d:s}. Moreover by the characterization of Schottky
groups acting on the Riemann sphere given by Maskit
\cite{maskit1}, it is not hard to prove that every group  of
M\"obius transformations which is   Schottky  in the sense of
definition \ref{d:s},  is a    Schottky group in $PSL_2(\Bbb{C})$.

\item
 In \cite{nori} Nori   gave the following construction of the  higher-dimensional
 analogues of the classical Schottky groups: let $n = 2k +1,\,k > 1$ and $g \geq 1$.
 Choose $2g$ mutually disjoint projective subspaces $L_1,\ldots,L_{2g}$ of dimension $k$ in
 $\Bbb{P}^n_{\Bbb{C}}$ and $0<\alpha<\frac{1}{2}$.  For every integer
 $ 1  \leq j \leq g$
choose a basis of $\Bbb{C}^{n+1}$ so that $L_j =[ \{z_0,\ldots,
z_k = 0\}-\{0\}]_n$ and $ L_{g+j} = [\{z_{k+1},\ldots, z_n =
0\}-\{0\}]_n.$ Define $\phi_j : \Bbb{P}^n_{\Bbb{{C}}} \rightarrow
\Bbb{R}$ by the formula $\phi_j [z_0,\ldots , z_n] =\frac {\mid
z_0\mid^2  + \ldots + \mid z_k\mid^2} {\mid z_{0}\mid^2  + \ldots
+ \mid z_{n}\mid^2} $ and consider the open neighborhoods $V_j =
\{x\in   \Bbb{P}^n_{\Bbb{C}} : \phi_j(x) < \alpha\}$ and $ V_{g+j}
= \{x\in\Bbb{P}^n_{\Bbb{C}} : \phi_j(x) > \alpha\}$ of $L_j$ and
$L_{g+j}$ respectively.  Consider the  automorphism  $\gamma_j$ of
$ \Bbb{P}^n_{\Bbb{C}}$ given by $ \gamma_j[z_0,\ldots, z_n] =
[\lambda z_0,\ldots,\lambda z_k, z_{k+1},\ldots, z_n]$ where
$\lambda\in \Bbb{C}$ and $\mid \lambda\mid =\frac{1}{\alpha}-1$.
Then $\gamma_j(V_j) = \Bbb{P}^n_{\Bbb{C}}- \overline{V_{g+j}}$.
Moreover  for all $\alpha$ small the group $\Gamma$ generated by
$\gamma_1,\ldots, \gamma_g$ is a Schottky group.

\item Let $L = \{L_1,\ldots,   L_g\}$, $g > 1$, be a set of $g$ projective
subspaces of dimension $n$ of $\Bbb{P}^{2n+1}_{\Bbb{C}}$,   all of them pairwise disjoint.
In \cite{sva}   it is shown that:
\begin{enumerate}

\item There  exists a set $\{V_1,\ldots, V_g\}$ of pairwise
disjoint open sets  of $\Bbb{P}^{2n+1}_{\Bbb{C}}$ such that
$\Bbb{P}^{2n+1}_{\Bbb{C}}-\partial (V_i)$ has 2 connected
components for each $1\leq j\leq g$,  $L_j$ is contained in  $V_j$
and the closures of the $g$ open sets are pairwise disjoint.

\item There are  involutions $T_1,\ldots,  T_g$ of
$\Bbb{P}^{2n+1}_{\Bbb{C}}$, such that each $T_j$, $j =
1,\ldots,g$, interchanges the two connected components of
$\Bbb{P}^{2n+1}_{\Bbb{C}}-\partial(V_j)$ and the boundary
$\partial(V_j)$ is an invariant set.

\item  Let $\Gamma$ be the group generated by $T_1,\ldots, T_g$
and let $\tilde \Gamma\cup\{id\}$ be the subgroup consisting  of
elements of $\Gamma$ which can be written as reduced words of even
length in the generators  (recall that $
 w=z_n^{\varepsilon_n}\cdots z_2^{\varepsilon_2}z_1^{\varepsilon_1}\in\Gamma$  is a
 reduced word of length $n$ if $ z_\ell\in \{T_1,\ldots T_g\}$;
 $ \varepsilon_\ell\in\{-1,+1\}$ and  if $ z_j=z_{j+1}$ then
 $ \varepsilon_j=\varepsilon_{j+1}$). For $g>2$ it is verified that  $\tilde \Gamma$ is a
 Schottky group in the sense of definition \ref{d:s}.
\end{enumerate}
\end{enumerate}
\end{examples}

We prove:

\begin{theorem} \label{t:schottky} If  $\Gamma \leq PSL_{2n+1}(\Bbb{C})$ is a discrete
subgroup, then $\Gamma$ cannot be a Schottky group acting on $\Bbb{P}^ {2n}_{\Bbb{C}}$.
\end{theorem}

\subsection{Basic Properties of Schottky Groups}
\label{sec:2}

\begin{definition}
For a subgroup  $\Gamma\leq PSL_n(\Bbb{C})$ satisfying
definition \ref{d:s}  we define:
\begin{enumerate}
\item
$F(\Gamma)=\Bbb{P}^n_{\Bbb{C}}-(\bigcup_{\gamma\in Gen(\Gamma)}R^*_\gamma\cup S^*_\gamma )$.
\item $ \Omega(\Gamma)=\bigcup_{\gamma\in \Gamma}\gamma (F(\Gamma))$.
\end{enumerate}
\end{definition}

\begin{example} If $\Gamma\leq PSL_{2n}(\Bbb{C})$ is any of the groups of the  example
\ref{e:grou}, then $\Bbb{P}^{2n+1}_{\Bbb{C}}-\Omega(\Gamma)$ is homeomorphic to
$\Bbb{P}^{n}_{\Bbb{C}}\times \mathcal{C}$, where $\mathcal{C}$ is a Cantor set, see
\cite{ maskit, nori, sva}.
\end{example}

\begin{proposition} \label{t:sa} If $\Gamma$ is a Schottky group, then:
\begin{enumerate}
\item $\Gamma$ is a free group  generated by $Gen(\Gamma)$.
\item $ \Omega(\Gamma)/\Gamma$ is a compact complex  $n$-manifold  and  $Int(F(\Gamma))$
is a fundamental domain for the action of $\Gamma$.
\end{enumerate}
\end{proposition}

Before we prove this result we  state a definition and prove a technical lemma.

\begin{definition}
Let $\Gamma\leq PSL_n(\Bbb{C})$ be a  subgroup. For an infinite
subset  $H\subset \Gamma$ and a  non-empty,  $\Gamma$-invariant
open set $\Omega\subset \Bbb{P}^n_{\Bbb{C}}$, we define
$Ac(H,\Omega)$ to be the closure   of the set of cluster points of
$HK$, where $K$ runs over all the compact subsets of $\Omega$.
Recall that  $p$ is a cluster point of $HK$ if there is a sequence
$(g_n)_{n\in\Bbb{N}}\subset H$ of different elements  and
$(x_n)_{n\in N}\subset K$ such that $g_n(x_n) \xymatrix{ \ar[r]_{n
\rightarrow  \infty}&} p$.
\end{definition}

\begin{lemma}  \label{l:lim}
For a subgroup  $\Gamma\leq PSL_n(\Bbb{C})$ satisfying definition
\ref{d:s} one has:
\begin{enumerate}
\item \label{i:lim1} \label{l:limp1} For each  reduced word  $
w=z_n^{\varepsilon_n}\cdots
z_2^{\varepsilon_2}z_1^{\varepsilon_1}\in\Gamma$  one has:
\begin{enumerate}
\item  If $\varepsilon_n =1$ then $w(Int(F(\Gamma)))\subset S^*_{z_n}$.
\item  If $\varepsilon_n =-1$ then $w(Int(F(\Gamma))\subset R^*_{z_n}$.
\end{enumerate}
\item \label{i:lim2} Let  $\gamma\in Gen(\Gamma)$. Then
$R(\gamma)=\bigcap_{k\in\Bbb {N}\cup\{0\}} \gamma^{-k}(R^*_\gamma)$ and
$ S(\gamma)=\bigcap_{k\in \Bbb{N}\cup \{0\}} \gamma^{k}(S^*_\gamma)$ are  closed disjoint
sets contained in $\Bbb{P}^{n}_{\Bbb{C}}-\Omega(\Gamma)$.

\item \label{i:lim3} Let $F_k=\{\gamma(f):f\in F(\Gamma)\textrm{ and } \gamma\in \Gamma  \textrm{ is  a reduced  word of Length at most } k \}$. Then
$F(\Gamma)\subset F_1(\Gamma)\subset\ldots \subset F_k(\Gamma)\subset \ldots$ and
$$\Omega(\Gamma) = \bigcup_{k\in \Bbb{N}\cup\{0\}}Int(F_k(\Gamma)). $$

\item \label{l:limp4} For each $\gamma\in Gen(\Gamma)$ one has
that $\emptyset\neq Ac(\{\gamma^n\}_{n\in
\Bbb{N}},\Omega(\Gamma))\subset S(\gamma)$ and $\emptyset\neq
Ac(\{\gamma^{-n}\}_{n\in\Bbb{N}},\Omega(\Gamma))\subset
R(\gamma)$.
\end{enumerate}
\end{lemma}

\begin{proof}

(\ref{i:lim1}) Let us  proceed by induction on the length of the  reduced words. Clearly
the case $k=1$ is done by the definition of Schottky group.  Now assume we have proven the
statement for $j=k$. Let $w=z_{k+1}^{\varepsilon_{k+1}}\cdots z_{1}^{\varepsilon_{1}}$ be
a reduced word and $x\in Int(F(\Gamma))$. By  the   inductive  hypothesis we deduce that
$z_{k+1}^{-\varepsilon_{k+1}}w(x)\in \Bbb{P}^{n}_{\Bbb{C}}-\overline{R_{z_{k+1}}^*}$ if
$\epsilon_{k+1}=1$ and $z_{k+1}^{-\varepsilon_{k+1}}w(x)\in
\Bbb{P}^{n}_{\Bbb{C}}-\overline{S_{z_{k+1}}^*}$ if $\epsilon_{k+1}=-1$. Now the  proof
follows by the definition of Schottky group.

(\ref{i:lim2})Let $\gamma\in Gen (\Gamma)$.  Since
$\gamma^m(\overline{S^*_\gamma})\subset \gamma^{m-1}(S^*_\gamma)$ we deduce that
$\bigcap_{m\in \Bbb{N}} \gamma^m(\overline S^*_\gamma)\subset \bigcap_{m\in \Bbb{N}}
\gamma^{m-1}(S^*_\gamma)= S(\gamma)$. To conclude observe  that:
$$ S(\gamma)=\bigcap_{m\in\Bbb{N}} \gamma^{m-1}(S^*_\gamma)\subset  \bigcap_{m\in \Bbb{N}}
\gamma^{m-1}(\overline S^*_\gamma)\subset \bigcap_{m\in \Bbb{N}}
\gamma^{m}(\overline
S^*_\gamma).$$

(\ref{i:lim3}) We will prove that that $F(\Gamma) \subset Int(F_1(\Gamma))$.  Let
$x\in \partial (F(\Gamma))$, then   there is  $\gamma_0\in Gen(\Gamma)$ such that
$x\in \partial  S^*_{\gamma_0}\cup \partial  R^*_{\gamma_0}$, for simplicity we will
assume that $x\in\partial  S^*_{\gamma_0}$. Define
$r_1=min  \{d(x,\gamma_0(\overline{S^*_\gamma })): \gamma\in Gen(\Gamma)\},\, r_2=min
\{d(x,\overline{R^*_\gamma}):\gamma\in Gen(\Gamma)\},\,
  r_3=min \{d(x,\gamma_0(\overline{R^*_\gamma})): \gamma\in Gen(\Gamma)-\{\gamma_0\} \},
  \,
r_4=min  \{d(x,\overline{S^*_\gamma}):\gamma\in
Gen(\Gamma)-\{\gamma_0\} \}$ and $r=min \{r_1,r_2,r_3,r_4\}$ (here
$d$ denotes the Fubini-Study metric). Clearly  $r>0$. Now,  let
$y\in B_{r/4}(x)\cap \overline{ S^*_{\gamma_0}}$ then by the
definition of $r$ we have that $y\in F(\Gamma)\cup \gamma
(F(\Gamma))$. If $y\in B_{r/2}(x)\cap \Bbb{P}^n_{\Bbb{C}}-
\overline{S^*_{\gamma_0}}$ then by definition of $r$ we deduce
$y\in F(\Gamma)$. In other words, we have shown  $F(\Gamma)\subset
Int(F_1(\Gamma))$. Therefore:
$$F_k(\Gamma)\subset \{\gamma(f): \gamma \textrm{ is  a reduced  word of length at most } k \textrm{ and } Int(F_1(\Gamma) \}  \subset F_{k+1}(\Gamma)$$  {\it i.e.}, $F_k(\Gamma)\subset Int(F_{k+1}(\Gamma))$. To
conclude the proof observe that:  $$\Omega(\Gamma)=\bigcup_{k\in
\Bbb{N}\cup \{0\}} F_k(\Gamma)\subset \bigcup_{k\in \Bbb{N}\cup
\{0\}} Int(F_{k+1}(\Gamma)) \subset\bigcup_{k\in \Bbb{N}\cup
\{0\}} Int(F_k(\Gamma)) \,.$$

(\ref{l:limp4})  Let $K\subset \Omega(\Gamma)$ be a compact set   and $x$ a cluster
point  of
$\{\gamma^m(K)\}_{m\in \Bbb{N}}$. Then there is a subsequence $(n_m)_{m\in\Bbb{N}}\subset
(m)_{m\in \Bbb{N}}$   and  a sequence $(x_m)_{m\in\Bbb{N}}\subset K$ such that
$\gamma^{n_m}(x_m)\xymatrix{ \ar[r]_{m \rightarrow  \infty}&}x$. In case
$x\notin  S(\gamma)$ it is deduced that    there is $k_0 \in \Bbb{N} $ such that
$x\notin \gamma^{k_0}(\overline{S^*_{\gamma}}) $. Taking  $r=d(x,
\gamma^{k_0}(\overline{S^*_{\gamma}}))$ we have that:
\begin{equation} \label{chido}
B_{r/2}(x)\cap \gamma^{k_0}(\overline{S^*_{\gamma}})=\emptyset.
\end {equation}

 On the other hand, observe that since $K$ is compact, by  part (\ref{i:lim3}) of the
 present lemma   there is $l_0\in \Bbb{N}$ such that $K\subset F_{l_0}(\Gamma)$; also
 observe  that since $(n_m)_{m\in \Bbb{N} }$  is an strictly increasing sequence, there is
 $k_1\in \Bbb{N}$ such that  $n_m>l_0+1+k_0$ for $m>k_1$.  With these facts in mind we
 deduce
 $\gamma^{l_0+1}(K)\subset \overline{S_{\gamma}^*}$ and therefore:
 $$\gamma^{n_m}(x_m)\in \gamma^{n_m-l_0-1}(\overline {S^*_{\gamma}})\subset
 \gamma^{k_0}(\overline {S^*_{\gamma}}) \textrm{ for }m>k_1. $$
 Hence
 $x\in \gamma^{k_0}(\overline {S^*_{\gamma}})$, which contradicts \ref{chido}.
 Thus $\emptyset\neq Ac(\{\gamma^{n}\}_{n\in\Bbb{N}},\Omega(\Gamma))\subset
 S(\gamma)$. Observe that similar arguments prove also
 $\emptyset\neq Ac(\{\gamma^{-n}\}_{n\in\Bbb{N}},\Omega(\Gamma))\subset R(\gamma)$.
\end{proof}

\textit{ Proof of proposition  \ref{t:sa}.}

(1) Assume   there is   a reduced word $h$  with length $>0$ such
that $h=Id$. Now, let   $x\in Int(F(\Gamma)) $, then by part
(\ref{l:limp1}) of lemma \ref{l:lim}, $x=h(x)\in
\bigcup_{\gamma\in Gen(\Gamma)}(R^*_{\gamma}\cup S^*_{\gamma})$,
which contradicts the choice of $x$. Therefore $\Gamma$ is free.

(2) Let $K\subset  \Omega(\Gamma) $ be a compact set, then by part
(\ref{i:lim3}) of lemma \ref{l:lim}, there is $k\in\Bbb {N}$ such
that $K\subset  F_k(\Gamma)$. Assume   there is  a word $w$ with
length $\geq 2k+2$ such that $w(F_k(\Gamma))\cap F_k(\Gamma) \neq
\emptyset$. So  there are  $x_1,x_2\in F(\Gamma)$ and   words
$w_1,\, w_2$ of length at most $k$ such that
$x_1=w_1^{-1}w^{-1}w_2x_2$. On the other hand  $w_1^{-1}w^{-1}w_2$
is a word with length $\geq 2$. By (\ref{l:limp1}) of lemma
\ref{l:lim}, $x_1=w_1^{-1}w^{-1}(w_2(x_2))\in \bigcup_{g\in
Gen(\Gamma)}S^*_{j}\cup R^*_{j}$, but this contradicts the choice
of $x_1$. Therefore $\Gamma$ acts properly discontinuously and
freely on $ \Omega (\Gamma)$. \qed

\begin{remark}
All the results in this section remain valid if we change $\Bbb{P}^{n}_{\Bbb{C}}$ for
 $\Bbb{P}^n_{\Bbb{R}}$.
\end{remark}

\section{Dynamics of Projective Transformations}

\begin{lemma} \label{l:ein}
Let $V$ be a $\Bbb{C}$-linear space with $dim_{\Bbb{C}}(V)=n$, $T:V\rightarrow V$  an
invertible  linear transformation and   $\lambda\in \Bbb{C}$ such that
$\mid \alpha\mid <\mid \lambda\mid $ for every eigenvalue $\alpha$ of $T$.
For every $l\in \Bbb{N}$ we have uniform convergence
$
\lambda^{-m}
\left(
\begin{array}{l}
m\\
l
\end{array}
\right )T^m \xymatrix{ \ar[r]_{m \rightarrow  \infty}&} 0 $ on
compact subsets  of $V$.
\end{lemma}

Here $ \left (
\begin{array}{l}
m\\
l
\end{array}
\right ) $ denotes  the number of sets with $l$ elements from a
set with  $m$ elements.

\begin{proof} Decomposing $T$ into one or more Jordan blocks according  to Jordan's
Normal Form Theorem we reduce the problem  to the case where
 there is $0<\mid \lambda\mid <1$ and  an ordered basis
$\beta=\{v_1,\ldots, v_n\}$, $n\geq 2$,  such that the matrix of
$T$ with respect to  $\beta$ (in symbols $[T]_{\beta}$) satisfies:

$$
[T]_{\beta}=
\left(
\begin{array}{lllll}
\lambda & 1      & 0      & \cdots & 0\\
0      & \lambda & 1      & \cdots & 0\\
\cdots & \cdots & \cdots & \cdots & \cdots \\
0      & 0      & 0      & \cdots & \lambda\\
\end{array}
\right).
$$
An inductive argument  shows thats for all $m>n$:
\begin{equation}  \label{e:jb}
[T^m]_{\beta}=
\left(
\begin{array}{ccccc}
\lambda^m & \left(\begin{array}{c}m\\1\end{array}\right)\lambda^{m-1} &
\left(\begin{array}{c}m\\2\end{array}\right) \lambda^{m-2}&\cdots &
\left(\begin{array}{c}m\\n-1\end{array}\right) \alpha^{m+1-n}\\
0    & \lambda^m & \left(\begin{array}{c}m\\1\end{array}\right)
\lambda^{m-1}&\cdots & \left(\begin{array}{c}m\\n-2\end{array}\right)
\lambda^{m+2-n}\\\cdots & \cdots & \cdots & \cdots & \cdots \\
0      & 0      & 0      & \cdots & \lambda^m\\
\end{array}
\right).
\end{equation}
For a compact subset  $K\subset V$  set $\sigma(K)=sup\{\sum_{j=1}^n\mid \alpha_j\mid:
\sum_{j=1}^n \alpha_jv_j\in K\}$. Let  $z\in K$,   $z=\sum_{j=1}^n\alpha_j v_j$, then
by equation \ref{e:jb} we deduce:

$$
\mid T^m(z)\mid \leq \sigma(K)max\{\mid v_j\mid :1\leq j\leq n\}
\sum_{j=1}^n\sum^{j-1}_{k=0}\left(\begin{array}{c}m\\k\end{array}\right )
\mid \alpha^{m-k} \mid,
$$
Hence it is sufficient to observe that:
$$
\left\vert
\left(
\begin{array}{l}
m\\
l
\end{array}
\right )
 \left ( \begin{array}{c}m\\k\end{array}\right)\alpha^{m-k} \right \vert
 \leq  m^{2max\{k,l\}} \mid\alpha\mid^{m-k}\xymatrix{ \ar[r]_{m \rightarrow  \infty}&} 0.
$$
\end{proof}

\begin{definition}
Let $V$ be a $\Bbb{C}$-linear space with $dim_{\Bbb{C}}(V)=n$ and  let $T:V\rightarrow V$
be a  $\Bbb{C}-$linear transformation. We define
$Eve(T)=\langle\langle\{v\in V: v \textrm{ is an eigenvector of } T\}\rangle \rangle.
$ Where   $\langle\langle \{v\in V: v \textrm{ is an eigenvector of } T\}\rangle\rangle$
will denote the linear subspace generated by the eigenvectos of $T$.
\end{definition}

\begin{lemma} \label{l:fc} Let $l,k\in \Bbb{N}\cup\{0\}$ with $l<k$. Then
$
\left (
\begin{array}{c}
m\\
l
\end{array}
\right )
\left (
\begin{array}{c}
m\\
k
\end{array}
\right )^{-1}
\xymatrix{ \ar[r]_{m \rightarrow  \infty}&}0.$
\end{lemma}

\begin{proof}

$\left(\begin{array}{c}m\\l\end{array}\right)\left(\begin{array}{c}m\\k
\end{array}\right)^{-1}=\prod_{j=l}^{k-1}\frac{j+1}{m-j}\leq \left(\frac{k}{m-l}
\right)^{k-l} \xymatrix{ \ar[r]_{m \rightarrow  \infty}&}0.
$
\end{proof}
\begin{lemma} \label{l:rt}
Let $V$ be a $\Bbb{C}$-linear space with $dim_{\Bbb{C}}(V)=n> 1$ and let
$T:V\rightarrow V$ be an invertible  linear transformation such that there  are
$ \lambda\in \Bbb{C}$, with $\mid \lambda\mid=1$, and an ordered  basis
$\beta=\{v_1,\ldots,v_n\} $ for which:
$$
[T]_\beta=
\left(
\begin{array}{lllll}
\lambda & 1      & 0      & \cdots & 0\\
0      & \lambda & 1      & \cdots & 0\\
\cdots & \cdots & \cdots & \cdots & \cdots \\
0      & 0      & 0      & \cdots & \lambda\\
\end{array}
\right)
$$ that is $[T]_{\beta}$ is a $n\times n$-Jordan block.
Then for every $v\in V-\{0\}$ there is a unique $k(v,T)\in \Bbb{N}\cup \{0\}$ such that
the set of cluster points of
$\left \{  \left ( \begin{array}{l}m\\k(v,T)\end{array}\right)^{-1}T^m(v)
\right\}_{m\in \Bbb{N}}$ lies in   $\langle\langle v_1\rangle\rangle-\{0\}$.
\end{lemma}

\begin{proof}  Let $z=\sum_{j=0}^n\alpha_j v_j$ and $k(z,T)=max\{1\leq j\leq n:
\alpha_j\neq 0\}-1$, then we have that:

$$
\left ( \begin{array}{l}m\\k(v,T)\end{array}\right)^{-1}T^m (z)
=
\sum_{j=1}^n
(
\sum^{n-j}_{k=0}
\left ( \begin{array}{l}m\\k\end{array}\right)
\left ( \begin{array}{l}m\\k(v,T)\end{array}\right)^{-1}\lambda^{m-k}\alpha_{k+j})v_j
$$
The result now follows  from lemma \ref{l:fc}.
\end{proof}

\begin{corollary} \label{l:rot}
Let $V$ be a $\Bbb{C}$-linear space with $dim_{\Bbb{C}}(V)=n$ and let
$T:V\rightarrow V$ be a  linear transformation such that there are
$ \alpha_1,\ldots,\alpha_n \in \Bbb{C}$, with $\mid \alpha_j\mid=1$ for each
$0\leq j\leq 0$, and an ordered  basis $\beta=\{v_1,\ldots,v_n\} $ for which
$T(\sum^{n}_{j=0}\beta_j v_j)=\sum^{n}_{j=0}\alpha_j\beta_j v_j$.
Then  $k(v,T)=0 $ is the unique positive integer for which the set of cluster points
of   $\left\{  \left ( \begin{array}{l}m\\k(v,T)\end{array}\right)^{-1}T^m(v)
 \right\}_{m\in \Bbb{N}}$, where $v\in V-\{0\}$,  lies on   $V-\{0\}$.
\end{corollary}

\begin{corollary} \label{c:uni}
Let $V$ be a $\Bbb{C}$-linear space with $dim_{\Bbb{C}}(V)=n$ and let  $T:V\rightarrow V$
be an invertible  linear transformation such that each of its eigenvalues is a unitary
complex number. Then for every $v\in V-\{0\}$ there is a unique $k(v,T)\in \Bbb{N}
\cup \{0\}$ for which the set of cluster points of $\left\{  \left (
\begin{array}{l}m\\k(v,T)\end{array}\right)^{-1}T^m(v)  \right\}_{m\in \Bbb{N}}$ lies in
 $Eve(T)-\{0\}$.
\end{corollary}
\begin{proof}
By the Jordan's Normal Form Theorem there are $k\in \Bbb{N}$; $V_1,\ldots,  V_k\subset V$
 linear subspaces and     $T_i:V_i\rightarrow V_i$, $1\leq i\leq k$  such that:
\begin{enumerate}
\item $\bigoplus_{j=1}^k V_j=V$.
\item For each  $1\leq i\leq k$,  $T_i$ is a non-zero $\Bbb{C}$-linear map whose
eigenvalues are unitary complex numbers.
\item $\bigoplus_{j=1}^kT_j=T$.
\item For each $1\leq i\leq k$, $T_i$ is either diagonalizable or $n_i=dim_{\Bbb{C}}>1$;
$V_i$ contains an ordered  basis $\beta_i$ for which $[T]_{\beta}$ is a
$n_i\times n_i$-Jordan block.
\end{enumerate}

Let $v\in V-\{0\}$ then there is a non empty finite set $W\subset \bigcup_{j=1}^k
V_j-\{0\}$ such that $v=\sum_{w\in  W} w$. Now,   take $i:W\rightarrow \Bbb{N}$ where
$i(w)$ is the unique element in $\{1,\ldots ,k \}$ such that $w\in V_{i(w)}$,
$k(v,T)=max\{k(w,T_{i(w)}):w\in W\}$, $W_1=\{w\in W: k(w,T_{i(w)})<k(v,T)\}$ and
$W_2=W-W_1$ then:

\begin{small}
\begin{equation} \label{e:suma}
\frac
{
T^m(v)
}
{
\left (
\begin{array}{l}
m\\
k(v,T)
\end{array}
\right )
}
=
\sum_{w\in W_1}
\frac
{
\left (
\begin{array}{l}
m\\
k(w,T_{i(w)})
\end{array}
\right )
}
{
\left (
\begin{array}{l}
m\\
k(v,T)
\end{array}
\right )
}
\frac
{
T^m_{i(w)}(w)
}
{
\left (
\begin{array}{l}
m\\
k(w,T_{i(w)})
\end{array}
\right )
}
+
\sum_{w\in W_2}
\frac
{
T^m_{i(w)}(w)
}
{
\left (
\begin{array}{l}
m\\
k(w,T_{i(w)})
\end{array}
\right )
}
\end{equation}
\end{small}

The result now follows from equation \ref{e:suma}, lemmas \ref{l:fc}, \ref{l:rt} and
corollary \ref{l:rot}.
\end{proof}

\begin{definition}
Let $\gamma\in PSL_n(\Bbb{C})$ be an element of infinite order and
let $\tilde \gamma$ be a lifting of  $\gamma$. Then we define:

\begin{enumerate}
\item $\mid Eva(\gamma)\mid=\{\mid \lambda\mid \in \Bbb{R}: \lambda \textrm{  is an
eigenvalue } \textrm{ of } \tilde\gamma\}$
\item $L_r(\gamma)=\langle[\{v\in \Bbb{C}^{n}: v \textrm{ is an  eigenvector of }
\tilde \gamma \textrm{ and } \mid \tilde\gamma(v)\mid=r\mid  v\mid\}]_n\rangle.$
\item  $L(\gamma)$ as the closure of accumulation points of  $\{\gamma^m(z)\}_{m\in
\Bbb{Z}}$  where  $z\in\Bbb{P}^n_{\Bbb{C}}$.
\end{enumerate}
Clearly parts 1 and 2 of  this definition do not depend on the
choice of $\tilde \gamma$.
\end{definition}

\begin{proposition} \label{l:lila}
 Let $\gamma\in PSL_ {n+1}(\Bbb{C})$ be an element of infinite order, then:
$$ L(\gamma)=
\bigcup_{r\in \mid Eva(\Gamma)\mid } L_r(\gamma).$$
\end{proposition}

\begin{proof} Since  $\bigcup_{r\in \mid Eva(\Gamma)\mid } L_r(\gamma)\subset L(\Gamma)$
is trivially verified, it is enough to check that
$L(\gamma)\subset\bigcup_{r\in \mid Eva(\Gamma)\mid }
L_r(\gamma)$. Let $\tilde \gamma$ be a lifting of $\gamma$, then
by the Jordan's Normal Form Theorem there are $k\in \Bbb{N}$;
$V_1,\ldots ,  V_k\subset \Bbb{C}^{n+1}$  linear subspaces;
$\gamma_i:V_i\rightarrow V_i$, $1\leq i\leq k$ and
$r_1,\ldots,r_k\in \Bbb{R}$ which satisfy:
\begin{enumerate}
\item $\bigoplus_{j=1}^k V_j=\Bbb{C}^{n+1}$.
\item For each  $1\leq i\leq k$,  $\gamma_i$ is a non-zero  $\Bbb{C}$-linear map whose
eigenvalues are unitary complex numbers.
\item $0<r_1<r_2<,\ldots,<r_k$.
\item $\bigoplus_{j=1}^kr_j\gamma_j=\tilde \gamma$.
\end{enumerate}
 In what follows $(\tilde\gamma,k,\{V_i\}_{i=1}^k,\{\gamma_i\}_{i=1}^k,\{r_i\}_{i=1}^k)$
 will be called a decomposition for $\gamma$. Now let $[v]_n\in\Bbb{P}^n_{\Bbb{C}}$, thus
  $v=\sum_{j=1}^kv_j$ where $v_j\in V_j$. Set $j_0=max\{1\leq j\leq k:v_j\neq 0\}$. One has:

\begin{equation} \label{e:elimina}
\left (\begin{array}{c}m\\ k(v_{j_0}, T_{j_0})
\end{array}\right )^{-1}\frac{\tilde{\gamma}^m(v)}{r_{j_0}^m}=\sum_{j=1}^k\left
(\begin{array}{c}m\\ k(v_{j_0}, T_{j_0})
\end{array}\right )^{-1}\frac{r_j^m\gamma_j^m(v_j)}{r_{j_0}^m}.
\end{equation}

By equation \ref{e:elimina}, lemma \ref{l:ein} and corollary  \ref{c:uni} we conclude
that
the set of cluster points of $\{ \gamma^{m}(v)\}_{m\in\Bbb {Z}}$ lies in
$[Eve(\gamma_{j_0})-0]_{n}=L_{r_{j_0}}(\gamma)$.
\end{proof}

\section{Proof of the Main Theorem}

\begin{lemma} \label{l:conlim}
Let $\Gamma\leq PSL_{2n+1}(\Bbb{C})$ be a group  and $\Omega $ a
non-empty, $\Gamma$-invariant   open set where $\Gamma$ acts
properly discontinuously and such that whenever $l$ is a
projective subspace contained in $\Bbb{P}^{2n}_{\Bbb{C}}-\Omega$
then $dim_{\Bbb{C}}(l)<n$. Then for every
 $\gamma\in \Gamma$ with infinite order there is a connected set $\mathcal{L}
 (\gamma)\subset Ac(\{\gamma\}_{m\in \Bbb{Z}},\Omega)\cup L(\gamma)$ such that
 $L(\gamma)\subset \mathcal{L}(\gamma)$.
\end{lemma}

\begin{proof}
  Let $\gamma\in \Gamma$ be an element with infinite order, and
   choose a
  decomposition
  $(\tilde\gamma,k,\{V_i\}_{i=1}^k,\{\gamma\}_{i=1}^k,\{r_i\}_{i=1}^k)$    for $\gamma$. Take
  $j_0=min\{1\leq j\leq k:\sum_{i=1}^jdim_{{C}}(V_i)\geq n+1\}$. From proposition
  \ref{l:lila} we can assume that $k\geq 2$. For the moment let us assume that
  $j_0\neq 1,\, k$. Observe that since $\sum_{i=1}^{j_0}dim_{{C}}(V_i)\geq n+1$
  we conclude that there is  $w=\sum_{i=1}^{j_0}w_j\in\bigoplus_{j=1}^{j_0}V_j$ non-zero,
  where $w_i\in V_i$, such that   $[w]_{2n}\in \Omega$ and since $\Omega$ is open we can
  assume that $w_{j_0}$ is
non-zero.  Now, let $z\in \bigoplus_{j>j_0}V_j-\{0\} $ then by lemma \ref{l:ein}
$$w_m(z)=\left [w+\left (\begin{array}{l}m\\k(w,\gamma_{j_0})
\end{array}\right )\sum_{j>j_0}\left (\frac{r_{j_0}}{r_j}\right)^{m}\gamma_j^{-m}(z_j)
\right]_{2n}\xymatrix{ \ar[r]_{m \rightarrow  \infty}&}[w]_{2n}$$
thus for $m(z)$ large $(w_m(z))_{m\geq m(z)}\subset\Omega$. On the
other hand, by corollary \ref{c:uni} there is an strictly
increasing  sequence $(n_m)_{m\in {N}}\subset {N}$ and $w_0\in
Eve(\gamma_{j_0})  -\{0\}$ such that:

$$
\left (
\begin{array}{l}
n_m\\
k(w_{j_0},\gamma_{j_0})
\end{array}
\right )^{-1}
\gamma^{n_m}_{j_0}(w_0)
\xymatrix{ \ar[r]_{m \rightarrow  \infty}&}
w_0.
$$
From here and lemma \ref{l:ein}  we deduce that:
$$
\gamma^{n_m}(w_{n_m})=
[
\left (
\begin{array}{l}
n_m\\
k(w_{j_0},\gamma_{j_0})
\end{array}
\right )^{-1} \sum_{j\leq j_0}(\frac{r_j}{r_{j_0}})^{n_m}\gamma_j^{n_m}(w_j)+z
]_{2n}
\xymatrix{ \ar[r]_{m \rightarrow  \infty}&}
[w_0+z]_{2n}.
$$

From here it follows that:

$$\bigcup_{j>j_0}L_{r_j}(\gamma)\subset \langle [w_0]_{2n},
[\bigoplus _{j>j_0}V_j-\{0\}]_{2n}\rangle\subset Ac(\{\gamma^{m} \}_{m\in \Bbb{Z}},
\Omega)\cup L(\gamma).$$ To conclude consider the following observations:

\begin{enumerate}
\item [Obs. 1] Observe that in  the previous  argument, the
assumption $j\neq k$ is not crucial, so for the case $j=1$ it is
verified  that there is $w_1\in  L(r_1)$ such that
$$\bigcup_{j>1}L(r_j)\subset \langle w_1,[\bigoplus _{j>1}V_i-\{0\}]_{2n}\rangle\subset
Ac(\{\gamma^{m}\}_{m\in \Bbb{Z}},\Omega)\cup L(\gamma)$$ thus in
case $j=1$ it is enough to take
 $$\mathcal{L}(\gamma)= \langle w_1,[\bigoplus _{j>1}V_j-\{0\}]_{2n}\rangle
 \cup L_{r_1}(\gamma).$$

\item[Obs. 2] \label{o:2} Applying the same argument to
$\gamma^{-1}$ in the case $j_0\neq 1,\,k $, it is deduced that
there is $v\in L(r_{j_0})$ such that:
 $$\bigcup_{j<j_0}L(r_j)\subset \langle v,[\bigoplus _{j<j_0}V_i-\{0\}]_{2n}\rangle
 \subset
Ac(\{\gamma^m\}_{m\in \Bbb{Z}},\Omega).$$ Therefore in this case it is enough to take
$$\mathcal{L}(\gamma)= <v,[\bigoplus _{j<j_0}V_j-\{0\}]_{2n}>\cup <[w_0]_{2n},
[\bigoplus _{j>j_0}V_j-\{0\}]_{2n}>\cup L_{r_{j_0}}(\gamma).$$

\item[Obs. 3]   To obtain the result  in the case $j=k$ it is
enough to apply the same argument  used in Obs. 1 to
$\gamma^{-1}$.
\end{enumerate}
\end{proof}

\begin{lemma} \label{c:lines}
If $\Gamma\leq PSL_{2n+1}(\Bbb{C})$ is a  Schottky group then
$\Bbb{P}^{2n}_{\Bbb{C}}-\Omega(\Gamma)$  does not contain   a
projective subspace $\mathcal{V}$ with
$dim_{\Bbb{C}}(\mathcal{V})\geq n$.
\end{lemma}

\begin{proof} If  $\mathcal{V}\subset \Bbb{P}^{2n}_{\Bbb{C}}-\Omega(\Gamma)$ is a
projective subspace with  $dim_{\Bbb{C}}(\mathcal{V})\geq n$, then:
\[
\mathcal{V}\subset \Bbb{P}^{2n}_{\Bbb{C}}-\Omega(\Gamma)=\Bbb{P}^{2n}_{\Bbb{C}}-
\bigcup_{\gamma\in \Gamma}\gamma(F(\Gamma))\subset \Bbb{P}^{2n}_{\Bbb{C}}-
 F(\Gamma)=\bigcup_{g\in Gen(\Gamma)}R^*_{\gamma}\cup S^*_{\gamma}.
\]

Since  $\mathcal{V}$ is connected and
$(\mathcal{V}\cap\cup_{\gamma\in Gen(\Gamma)} R^*_\gamma,
\mathcal{V}\cap\cup_{\gamma\in Gen(\Gamma)} S^*_\gamma)$ is a
disconnection for $\mathcal{V}$ we deduce that
$\mathcal{V}\subset \cup_{\gamma\in
 Gen( \Gamma)R^*_\gamma}$ or $\mathcal{V}\subset \cup_{\gamma\in Gen( \Gamma)S^*_\gamma}$.
 Moreover by an inductive  argument we deduce that there is $ \gamma_0\in Gen(\Gamma)$
  such that $\mathcal{V}\subset S^*_{\gamma_0}$ or $\mathcal{V}\subset R^*_{\gamma_0}$.
  For simplicity let us assume that    $\mathcal{V}\subset S^*_{\gamma_0}$. Taking
  $\sigma\in Gen(\Gamma)-\{\gamma_0\}$ we have:
\begin{equation} \label{e:inlin}
\sigma^{-1}(\mathcal{V})\subset \sigma^{-1}( S^*_{\gamma_0}) \subset  \sigma^{-1}
(\Bbb{P}^{2n}_{\Bbb{C}}- \overline S^*_{\sigma} )=R^*_\sigma.
\end{equation}
 Observe that $\mathcal{V}$ and $\sigma^{-1}\mathcal{V}$ are projective  subspaces with
$dim_{\Bbb{C}}(\mathcal{V})+dim_{\Bbb{C}}(\sigma^{-1}\mathcal{V})\geq 2n$ then
$\mathcal{V} \cap \sigma^{-1} (\mathcal{V})\neq \emptyset$. However, this is a
contradiction since by equation \ref{e:inlin} we have that $\mathcal{V}\cap
\sigma^{-1}\mathcal{V}\subset R^*_\sigma \cap S^*_{\gamma_0} =  \emptyset$
\end{proof}

\textit {Proof of Theorem \ref{t:schottky}.} Assume that there is  a group
$\Gamma\leq PSL_{2n+1}(\Bbb{C})$ which is a Schottky group and let $\gamma\in Gen(\Gamma)$.
 By lemma \ref{l:conlim} there is  a connected  set $\mathcal{L}(\gamma)$ such that
 $L(\gamma)\subset\mathcal{L}(\gamma)\subset Ac(\{\gamma^m\}_{m\in \Bbb{N}},
 \Omega(\Gamma))$. On the other hand  by (\ref{l:limp4}) of lemma \ref{l:lim} we
 have $Ac(\{\gamma^m\}_{m\in \Bbb{N}},\Omega(\Gamma))\subset S( \gamma)\cup R(\gamma)$.
 Since
$(R(\gamma)\cap \mathcal{L}(\gamma),\,S(\gamma)\cap
\mathcal{L}(\gamma))$ is a disconnection for $\mathcal{L}(\gamma)$ we deduce   $\mathcal{L}(\gamma)\subset
R(\gamma)$
 or $\mathcal{L}(\gamma)\subset S(\gamma)$. This implies  $L(\gamma)\cap
 S(\gamma)=\emptyset$
or $L(\gamma)\cap R(\gamma)=\emptyset$. However this contradicts (\ref{l:limp4}) of
lemma \ref{l:lim}. Therefore  $\Gamma$ cannot be a Schottky group. \qed

\begin{remark}
\begin{enumerate}
\item If in   definition \ref{d:s} we allow  that $R^*_j=S^*_j$
and  $\gamma_j^2=Id$ for $1\leq j\leq g$,   the resulting group is
a type of Complex Kleinian Group (see \cite{sva}), and by means of
theorem \ref{t:schottky} is not hard to see that for $g\geq 3$
this type of groups cannot be realized as subgroups of
$PSL_{2n+1}(\Bbb{C})$. \item Theorem \ref{t:schottky} remains
valid if we change $\Bbb{C}$ by  $\Bbb{R}$.
\end{enumerate}
\end {remark}

\begin{center}
{\sc Acknowledgements}
\end{center}
The author would like to thank to the  Instituto de Matem\'aticas Unidad Cuernavaca de la  
Universidad Nacional Autonoma de M\'exico for his hospitality during the writing of this paper, to A. Verjovsky for  suggesting the problem and  his valuable discussions,  to J. Seade for his the advice and valuable discussions during the writing of this paper.



\end{document}